\magnification 1180%\magstep1
\input epsf.sty
\input  amssym.tex
%\headline{\hfill November 3, 2007}
\overfullrule0pt
%\raggedbottom
  
    %PS-Palatino-Bold
  % rm for letterhead
%\font\rm=cmr10 at 10pt  % rm for text

  % small rm for letterhead 
%\font\rm=PS-Bookman-Light at 11pt   % rm for text

%\font\smtt=tt10 at 8 pt

\def\sqr#1#2{{\vcenter{\hrule height.#2pt              %qed
     \hbox{\vrule width.#2pt height#1pt\kern#1pt
     \vrule width.#2pt}
     \hrule height.#2pt}}}
\def\square{\mathchoice\sqr{5.5}4\sqr{5.0}4\sqr{4.8}3\sqr{4.8}3}
\def\qed{\hskip4pt plus1fill\ $\square$\par\medbreak}

\def\P{{\bf P}}
\def\C{{\bf C}}

\vskip.3in
\centerline{\bf Degree Complexity of a Family of Birational Maps}

\bigskip
\centerline{ Eric Bedford\footnote{*}{Research supported in part by the NSF}, Kyounghee Kim,  Truong Trung Tuyen,}

\centerline{Nina Abarenkova\footnote{\dag}{Research supported in part by a Russian Academy of Sciences/CNRS program.}, and Jean-Marie Maillard}

\bigskip

\noindent {\bf \S0.  Introduction.  }  Birational mappings on the space of $q\times q$ matrices have been found to arise as natural symmetries in lattice statistical mechanics.  One such map gives rise to a family $k_{a,b}$ of birational maps of the plane (see [BMR1,2], [A2]).  Dynamical properties of this family have been studied in a number of works ([A1--8], BD[2], [BMR], [BM]). Recall the quantity
$$\delta(k):=\lim_{n\to\infty}({\rm deg}(k^n))^{1\over n},$$
which is the exponential rate of growth of the iterates of $k$.  This is variously known as the degree complexity, the dynamical degree, or the algebraic entropy  of $k$.  When $b\ne0$ and $a$ is generic, $\delta(k_{a,b})$ is the largest root of the polynomial $x^3-x^2-2x-1$.  When $b=0$ and $a$ is generic, $\delta(k_{a,0})$  is the largest root of $x^2-x-1$.   The form of a map can change radically under birational equivalence: a simpler form for $k_{a,0}$ which was obtained in [BHM] made it more accessible to detailed analysis (see [BD1,3]).

A basic property is that $k$ is reversible in the sense that $k=\jmath\circ\iota$ is a composition of two involutions.  Here we give (a  birationally equivalent version of) $k$  as a composition of involutions in a new way.  This shows how $k_{a,b}$ fits naturally into a larger family of maps.  Namely, for any polynomial $F$, we define the involutions
$$\jmath_F(x,y)=(-x+F(y),y),\ \ \ \iota(x,y) = \left(1-x-{x-1\over y},-y-1-{y\over x-1}\right),$$
and the family of birational maps is given by $k_F=\jmath_F\circ \iota$.   When $F$ is constant, the family $k_F$ is birationally equivalent to $k_{a,0}$, and when $F$ is linear, $k_F$ is equivalent to $k_{a,b}$.  In this paper we determine the structure and degree complexity for the maps $k_F$:
\proclaim Theorem 1.  Let $n$ denote the degree of $F$.  If $n$ is even, then for generic parameters $\delta(k_F)$ is the largest root of the polynomial $x^2-(n+1)x-1$.  If $n$ is odd, then for generic parameters $\delta(k_F)$ is the largest root of $x^3-nx^2-(n+1)x-1$.

When $k_F$ is not generic, the growth rate $\delta(k_F)$ decreases (i.e.\ $F\mapsto\delta(k_F)$ is lower semicontinuous in the Zariski topology).  One of the interesting things about the family is to know which parameters are not generic as well as the corresponding values of $\delta(k_F)$  is decreased.  The exceptional values of $a$ for the family $k_{a,0}$, as well as the corresponding values of $\delta(k_{a,0})$, were found  by Diller and Favre [DF].  Similarly, the exceptional values of $(a,b)$ are given in [BD2].  Here we look at the maximally exceptional parameters for the case where $F$ is cubic. These are the cubic maps with the slowest degree growth and give a 2 complex parameter family of maps which are (equivalent to) automorphisms: 

\proclaim Theorem 2.  If $F(y)= ay^3+ay^2+by+2$, $a\ne0$, then $k_F$ is an automorphism of a compact, complex surface ${\cal Z}$.  Further, the degrees of $k^n_F$ grow quadratically, and  $k_F$ is integrable.

We will analyze the family $k_F$ by inspecting the blowing-up and blowing-down behavior.  That is, there are exceptional curves, which are mapped to points; and there are points of indeterminacy, which are blown up to curves.  As was noted by Forn\ae ss and Sibony [FS], if there is an exceptional curve whose orbit lands on a point of indeterminacy, then the degree is not multiplicative:  $({\rm deg}(k_F))^n\ne {\rm deg}(k_F^n)$.    The approach we use here is to replace the original domain ${\bf P}^2$ by a new manifold ${\cal X}$.  That is, we find a birational map $\varphi:{\cal X}\to{\bf P}^2$, and we consider the new birational map $\tilde k=\varphi\circ k_F\circ \varphi^{-1}$.   There is a well defined map ${\tilde k}^*:Pic({\cal X})\to Pic({\cal X})$, and the point is to choose ${\cal X}$ so that the induced map $\tilde k$ satisfies $({\tilde k}^*)^n =({\tilde k}^n)^*$.  By the birational invariance of $\delta$ (see [BV] and [DF]) we conclude that $\delta(k_F)$ is the spectral radius of ${\tilde k}^*$.  This method has also been used by Takenawa [T1--3].  The  general existence of such  a map $\tilde k$ when $\delta(k)>1$ was shown in [DF].  We comment that the construction of ${\cal X}$ and $\tilde k$ can yield further information about the dynamics of $k$ (see, for instance, [BK] and [BD2]).

\bigskip 

\noindent{\bf\S1. The maps. } Let us set $F(z) = \sum_{j=0}^n a_j z^j$ with $a_n \ne 0$. The map $k=\jmath_F\circ \iota$ is  the composition of the two involutions defined above.  The map $k = [k_0:k_1:k_2]$ is given in homogeneous coordinates as
$$\eqalign{&k_0 =(x_0x_1-x_0^2)^n x_2 \cr& k_1= x_0^{n-1}(x_1-x_0)^{n+1} (x_2+x_0) +x_2 \sum_{j=0}^n a_j (x_0 x_1 -x_0^2)^{n-j} (x_2^2-x_0 x_1-x_1 x_2)^j \cr&k_2= x_2 (x_0x_1-x^2_0)^{n-1}(x_2^2-x_0x_1-x_1x_2).}\eqno(1.1)$$
Each coordinate function has degree $2n+1$, which means that ${\rm deg}(k) = 2 \cdot {\rm deg} (F) +1$. Since the jacobian of this map is $x_0^{3n-3} (x_0-x_1)^{3n-1} x_2^2 (x_0^2-x_0x_1-x_1x_2)$ we have four exceptional curves : 
$$C_1:=\{x_0=0\},\ C_2:=\{x_0=x_1\},\ C_3:=\{x_2=0\},\ C_4:=\{-x_0^2+x_0x_1+x_1x_2=0\}.$$
When $a_0 \ne 2$, the exceptional hypersurfaces are mapped as:
$$k \ :\   C_4 \mapsto [1:-1+a_0:0] \in C_3{\rm  \ \ \ \ and\ \ \ } C_1 \cup C_2\cup C_3 \mapsto e_1.\eqno{(1.2)}$$
The points of indeterminacy for $k$ are $$e_1:=[0:1:0],\ \ e_2 :=[0:0:1],\ \ {\rm and}\ \ e_{01}:=[1:1:0].$$
Figure 1.1 shows the relative position of the points of indeterminacy (dots with circles around them), exceptional curves, and the critical images (big dots).  The information that $C_1,C_2,C_3\to e_1$ is not drawn for lack of space.
\medskip
\epsfysize=1.5in
\centerline{ \epsfbox{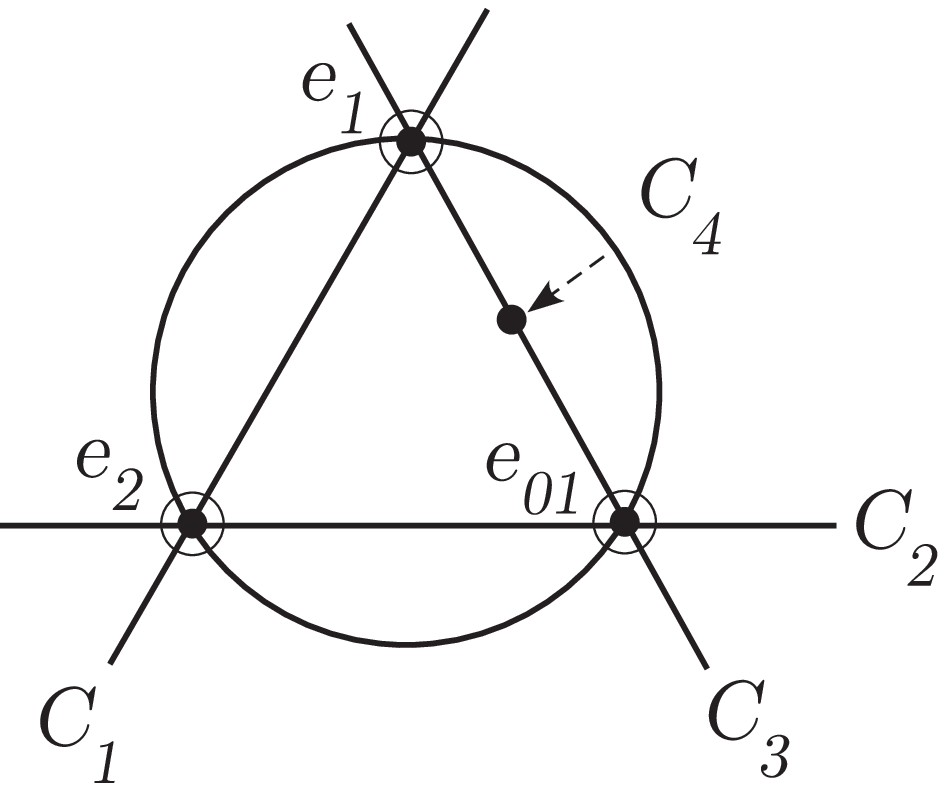}  }
\centerline{Figure 1.1.  Exceptional curves and points of indeterminacy. }

\medskip
The sort of singularity that will be the most difficult to deal with arises from the exceptional curve $C_1\mapsto e_1\in C_1$.  In local coordinates near $e_1$, this looks like
$$k[t:1:y] = \left[{t^n+\cdots\over a_n (-y)^n +\cdots}:1:{t^{n-1}+\cdots\over a_n (-y)^{n-1}+\cdots}\right].\eqno(1.3)$$
For this, we will perform the iterated blowups described in \S2.

The inverse map $k^{-1}=[k_0^{-1}:k_1^{-1}:k_2^{-1}]$ is given as
$$\eqalign{&k_0^{-1}= x_0^n x_2(\check F - x_0^{n-1} (x_0+x_1))\cr
&k_1^{-1} =(x_0+x_2) \left(\sum_{j=0}^n a_j x_0 ^{n-j} x_2 ^j - x_0^{n-1} (x_0+x_1)\right)^2\cr 
& k_2^{-1}= x_0^{n-1} x_2 \left(x_0^{n-1} (x_0^2+x_0x_1+x_1 x_2) - (x_0+x_2)\check F\right)} $$
where $\check F=x_0^nF(x_2/x_0)=\sum_{j=0}^n a_j x_0 ^{n-j} x_2 ^j$.
The jacobian for the inverse map is 
$$ x_0^{3n-3} x_2^2 (x_0^n +x_0^{n-1} x_1-\check F)^2 \left(x_0^{n+1} - (x_0+x_2) (x_0^n + x_0^{n-1} x_1 -\check F)\right)$$
The exceptional curves for $k^{-1}$ are $C_j'$, $1\le j\le 4$,  where 
$$\eqalign{&C_1'=C_1,\ \ C_2':=\{x_0^n +x_0^{n-1} x_1-\check F=0\},\ \ C_3'=C_3,  \cr  
&C_4':=\{x_0^{n+1} - (x_0+x_2) (x_0^n + x_0^{n-1} x_1 - \check F)=0\} .}$$
$$k^{-1} \ :\ C_1'\cup C_3' \mapsto e_1,\ \ C_2' \mapsto e_2,\ \ {\rm and\ \ }C_4' \mapsto e_{01} \in C_3',\ \ \eqno{(1.4)}$$

\noindent{\bf\S2. Blowups and local coordinate systems.}
In this section we discuss iterated blowups, and we explain the choices of local coordinates which will be useful in the sequel.  Let $\pi : X \to \C^2$ denote the complex manifold obtained by blowing up the origin $e=(0,0)$; the space is given by
$$X = \{ ((t,y),[\xi: \eta]) \in \C^2 \times \P\  ;\  t \eta = y \xi \},$$
and $\pi$ is projection to ${\bf C}^2$.
Let $E:=\pi^{-1} (e)$ denote the exceptional fiber over the origin, and note that $\pi^{-1}$ is well defined over ${\bf C}^2-e$. The closure in $X$ of the $y$-axis ($\pi^{-1}(\{t=0\}-e)$)  corresponds to the hypersurface $\{\xi=0\}\subset X$.  On the complement $\{\xi\ne0\}$ set $u=t$ and $\eta= y/t$. Then $(u,\eta)$ defines a coordinate system  on $X\setminus \{t=0\}$, with a point being given by $((t,y),[1:y/t])=((u,u\eta),[1:\eta])$.  We will use the notation $(u,\eta)_L$.  On the set $t\ne0$, the coordinate projection $\pi$ is given in these coordinates as
$$\pi_L(u, \eta)_L = (u, u \eta)=(t,y) \in \C^2.\eqno(2.1)$$  

Figure 2.1 illustrates this blowup with emphasis on the relation between the point $e$ and the lines $t=0$ and $y=0$ which contain it.  The space $X$ is drawn twice to show two choices of coordinate system; the dashed lines show where each coordinate system fails to be defined.  The left hand copy of $X$ shows the $u,\eta$-coordinate system in the complement of $t=0$.  The right hand side shows a different choice of coordinate; we would choose this coordinate system to work in a neighborhood of the point $p_1:=E\cap\{t=0\}$.

In the $u,\eta$ coordinate system (on the upper left side of Figure 2.1), the $\eta$-axis $(u=0)$ represents the exceptional fiber $E\cong\P^1$. The line $\gamma_\eta = \{ (s, \eta)_L : s \in \C\}$  projects to the line $\{y=\eta t\}\subset\C^2$, and  $(0,\eta)_L=E\cap \gamma_\eta$. It follows that $E\cap \{ y = 0 \} = (0,0)_L$ in this coordinate system. 
\medskip
\epsfysize=2.9in
\centerline{ \epsfbox{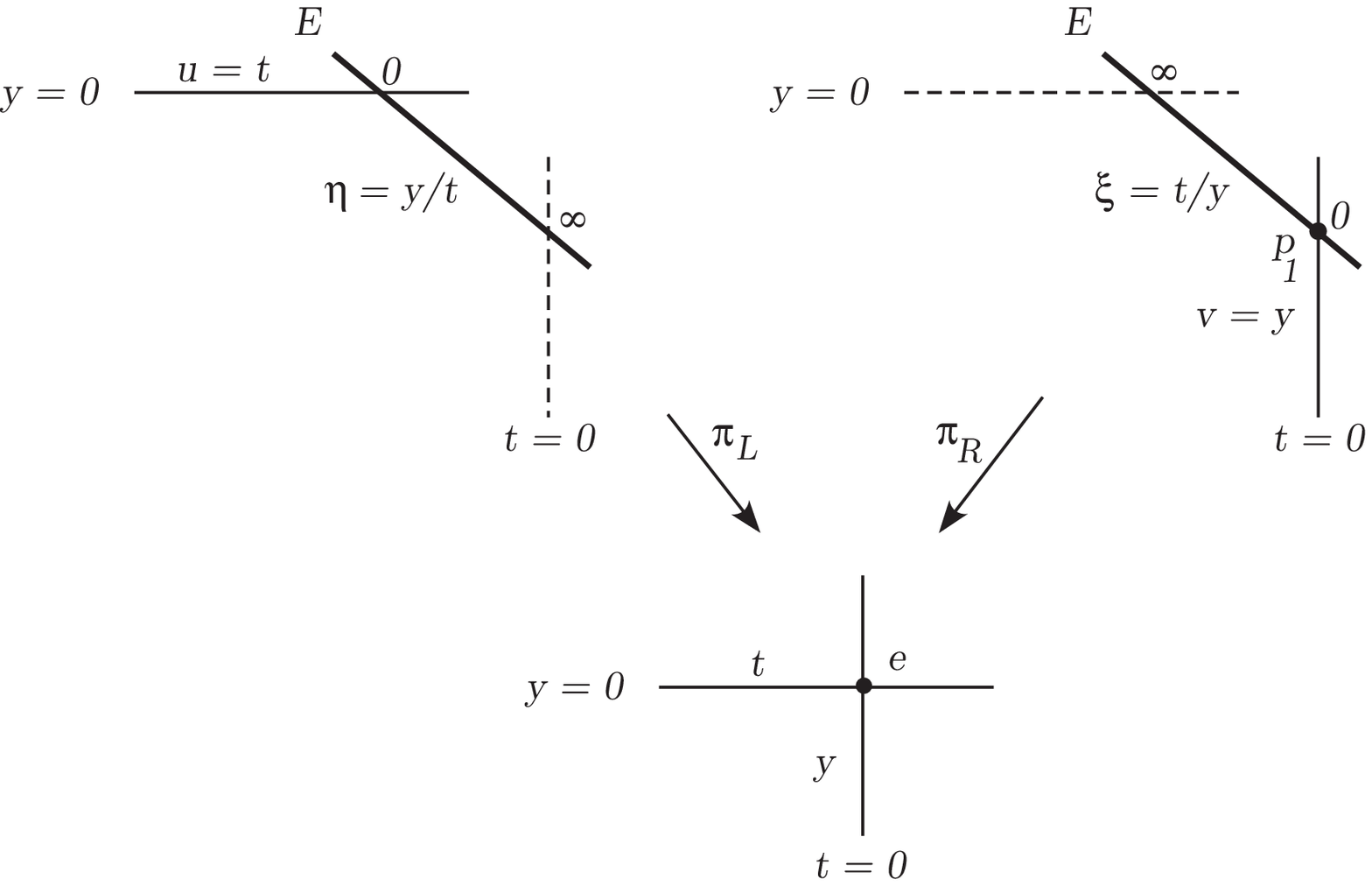}  }
\centerline{Figure 2.1.  Two choices of local coordinate systems. }
\medskip

\noindent On the upper right side of Figure 2.1, we define a $(\xi,v)$-coordinate system on the complement of $t$-axis ($y=0$) : 
$$\pi_R\ :\ (\xi,v)_R= ({t/ y},y) \to (v\xi, v) \in \C^2. \eqno{(2.2)}$$
The exceptional fiber $E$ is given by $\xi$-axis ($v=0$). 
Next we blow up $p_1=E\cap \{t=0\}= \{\xi=v=0\}=(0,0)_R$. Let $P_1$ denote the exceptional fiber over $p_1$. The choice of a local coordinate system depends on the center of next blowup. Suppose the third blowup center is an intersection of two exceptional fibers $p_2:=E \cap P_1$.  For this we are led to the $(u,\eta)$- coordinate system, as on the left side of Figure 2.1.   Thus we have a local coordinate system on the complement of $\{t=0\} \cup \{ y=0\}$;
$$ (u_1, \eta_1)_{1} = ({t/ y}, { y^2 / t} ) \to(u_1, u_1 \eta_1)_R\to (u_1^2 \eta_1, u_1 \eta_1) \in \C^2. \eqno{(2.3)}$$ 
This $(u_1, \eta_1)$-coordinate system is defined only off the axes $(t=0)\cup(y=0)$; the new exceptional fiber $P_1$ is given by the $\eta_1$-axis. 

\medskip
\epsfysize=1.5in
\centerline{ \epsfbox{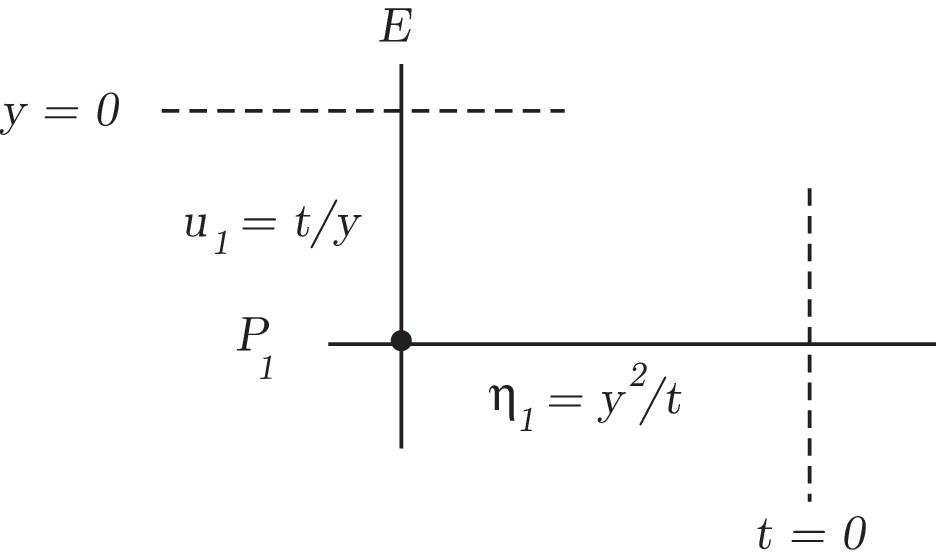}  }
\centerline{Figure 2.2.  Blowup of $p_1$ in $(u_1, \eta_1)$-coordinates. }
\medskip
Now we define a sequence of iterated blowups which will let us deal with the singularity (1.3).  We start with the blowup space $X$ as in Figure 2.2, and we continue inductively for  $2 \le j\le n$ by setting $p_j := E \cap P_{j-1}$ and letting $P_j$ be the exceptional fiber.   For each  $2 \le j\le n$, we use the left-hand coordinate system of Figure 2.1, which corresponds to (2.1).  Thus we have the coordinate projection $\pi_j:P_j\to{\bf C}^2$:
$$\pi_j:(u, \eta)_j \to(u^{j+1}\eta,u^j\eta)=(t,y)\in{\bf C}^2, \ \ \ \pi_j^{-1}(t,y)=(u,\eta)=(t/y,y^{j+1}/t^{j}). \eqno{(2.4)}$$ 
This coordinate system is defined off of $\{y=0\} \cup\{t=0\} \cup P_1 \cup \cdots \cup P_{j-1}$. A point $(0,\eta=c)_j\in P_j$ is the landing point of the curve $u\mapsto(u,c)_j$ as $u\to0$, which projects to the curve $u\mapsto (t(u)=u^{j+1}c,y(u)=u^{j}c)\in{\bf C}^2$.   In Figure 2.3, the exceptional fibers $P_j, \ 1 \le j \le n$ are drawn with their fiber coordinates $y^{j+1}/t^j$.

\medskip
\epsfysize=1.96in
\centerline{ \epsfbox{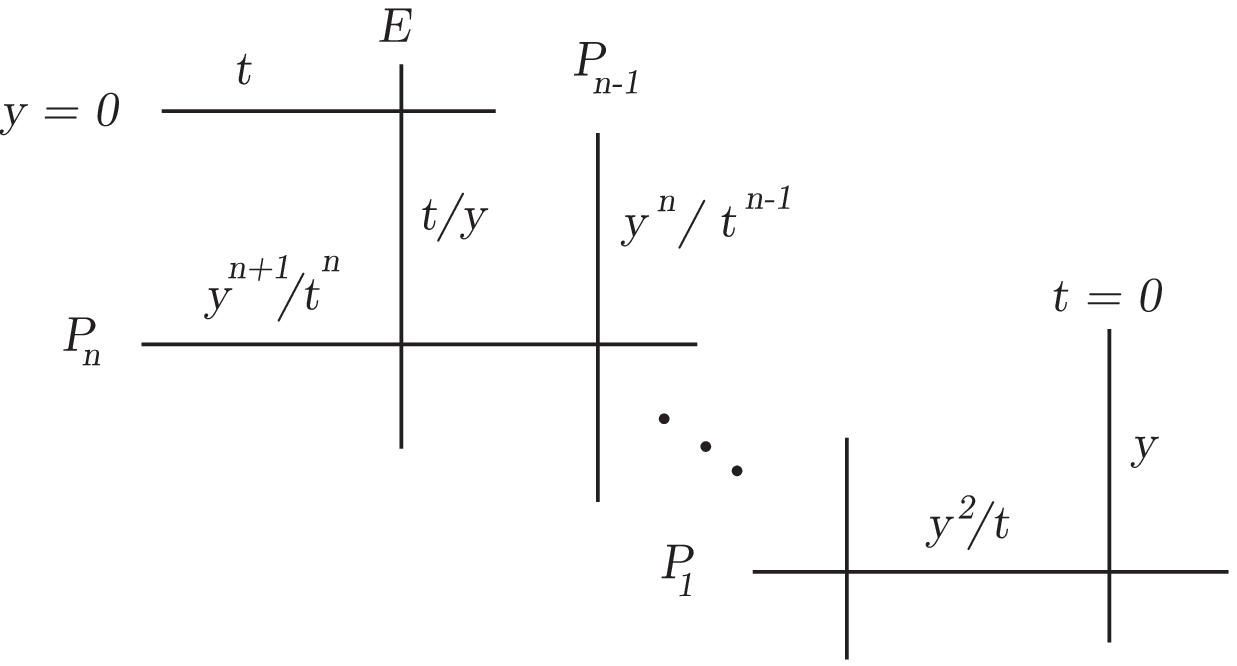}  }
\centerline{Figure 2.3.  $n$-th iterated blowup}
\medskip

\bigskip

\noindent{\bf\S3. Mappings with  ${ n=}$ even.}   We define  a complex manifold $\pi_{\cal X}: {\cal X} \to \P^2$ by blowing up points $e_1, q, p_1, \dots, p_{n-1}$  in the following order: \smallskip
\item{(i)} blow up $e_1=[0:1:0]$  and let $E_1$ denote the exceptional fiber over $e_1$,
\item{(ii)} blow up  $q:= E_1 \cap C_4$ and let $Q$ denote the exceptional fiber over $q$,
\item{(iii)} blow up $p_1 : = E_1 \cap C_1$ and let $P_1$ denote the exceptional fiber over $p_1$,
\item{(iv)} blow up $p_j :=E_1 \cap P_{j-1}$ with exceptional fiber $P_j$  for $2 \le j \le n-1$.
\smallskip
The iterated blow-up of $p_1,\dots,p_{n-1}$ is exactly the process described in \S2, so we will use the local coordinate systems defined there. That is, in a neighborhood of $Q$ we use a $(\xi_1,v_1)=(t^2/y,y/t)$ coordinate system. For $E_1$ and $P_j, 1 \le j\le n-1$ we use local coordinate systems defined in (2.2--4). We use homogeneous coordinates by identifying a point $(t,y)\in \C^2$ with $[t:1:y]\in \P^2$.  Let $k_{\cal X}:{\cal X}\to {\cal X}$ denote the induced map on the complex manifold ${\cal X}$. In the next few lemmas, we will show that $k_{\cal X}$ maps the exceptional fibers as shown in Figure 3.1.

\medskip
\epsfysize=1.7in
\centerline{ \epsfbox{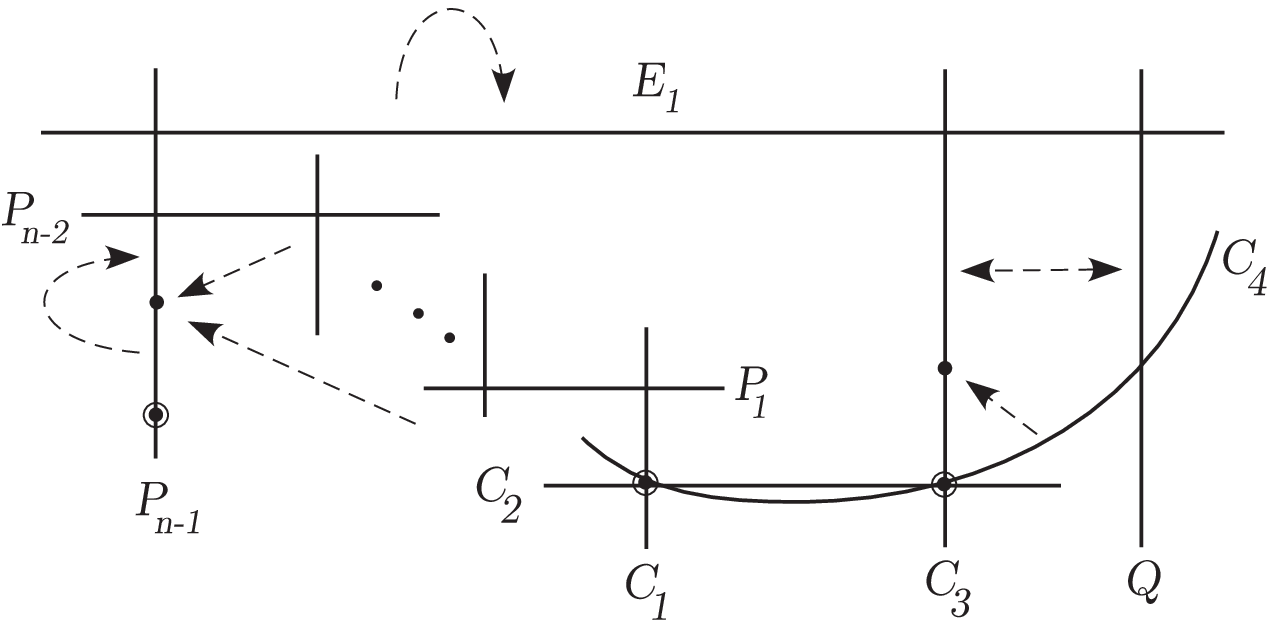}  }
\centerline{Figure 3.1.  The space ${\cal X}$  and the action of $k_{\cal X}$}
\medskip

\proclaim Lemma 3.1. Under the induced map $k_{\cal X}$, the blowup fibers $E_1$ and $P_{n-1}$ are mapped to themselves : $$\eqalign{k_{\cal X} \ \ :\ \  &E_1 \ni \xi\  \mapsto\  -\xi/(\xi+1) \in E_1\cr & P_{n-1} \ni \eta_{n-1} \ \mapsto \ \eta_{n-1}/(1+ a_n \eta_{n-1}) \in P_{n-1}.}\eqno{(3.1)}$$

\noindent{\it Proof.} First let us work on $E_1$.  We use the local coordinate system defined in (2.2), so a point in the exceptional fiber $E_1$ is $(\xi,0)_R$. To see the forward image of $E_1$ we consider a nearby point $(\xi, v)_R\to (v \xi, v)$ with small $v$ and we have $k_{\cal X} (\xi,0)_R = \lim_{v \to 0} k_{\cal X}(\xi, v)_R$.  By (1.1) we see that 
$$k \,[ v \xi:1:v] = [ v \xi + \cdots \,:\, 1+ \cdots\,:\, -v (\xi+1) +\cdots]$$ 
where we use $\cdots$ to indicate the higher order terms in $v$.  As in Figure 2.1, the coordinate of the landing point in $E_1$ is given by the ratio of $t$- and $y$-coordinates.  Thus we have $$k_{\cal X}|E_1 \ :\ \xi \mapsto  \lim_{v \to 0}k_0/k_2=\lim_{v \to 0} (v \xi + \cdots)/(-v (\xi+1) +\cdots) = - \xi/(\xi+1).$$

Now we determine the behavior of $k_{\cal X}$ on $P_{n-1}$.  A fiber point  $(0,\eta_{n-1})\in P_{n-1}$ is the landing point of the arc $u\mapsto(u,\eta_{n-1})$ as $u\to 0$.  To show that $k_{\cal X}$ maps $P_{n-1}$ to $P_{n-1}$, we need to evaluate:
$$\lim_{u\to0}k_{\cal X}(u,\eta_{n-1})=\lim_{u\to0}\pi_{n-1}^{-1}\circ k\circ\pi_{n-1}(u,\eta_{n-1}).$$
Using the formulas for $\pi_{n-1}$ and $\pi_{n-1}^{-1}$ in (2.4), we obtained the desired limit.  \qed

Now we may use similar calculations to show that $k_{\cal X}:P_j\to P_{n-1}$; we fix a point $(0,\eta_j)\in P_j$ and show the existence of the limit
$$\lim_{u\to 0}k_{\cal X}(u,\eta_j)=\lim_{u\to0}\pi_{n-1}^{-1}\circ k\circ \pi_j(u,\eta_j).$$
Doing this, we find that the line $C_1$ and all blowup fibers $P_j, j=1, \dots, n-2$ are all exceptional for both $k_{\cal X}$ and $k_{\cal X}^{-1}$. And $C_2$ is exceptional for $k_{\cal X}$:
$$\eqalign{&k_{\cal X}\ :\ C_1, C_2, P_1 , \cdots , P_{n-2} \mapsto 1/a_n \in P_{n-1} \cr & k_{\cal X}^{-1}\ :\ C_1 , P_1, \cdots , P_{n-2} \mapsto (-1)^{n-1}/a_n \in P_{n-1} } \eqno{(3.2)}$$
Combining (3.1--2) it is clear that the indeterminacy locus of $k_{\cal X}$ consists of three points 
$$e_2, \ \ e_{01},\ \  {\rm\ and\ \ } (-1)^{n-1}/a_n \in P_{n-1}.$$

\proclaim Lemma 3.2. If $n$ is even, then the orbits of the exceptional curves $C_1, C_2, P_1, \dots, P_{n-2}$ are disjoint from the indeterminacy locus. 

\noindent{\it Proof.} By Lemma 3.1,  the orbit of $1/a_n$ in $P_{n-1}$ is  $\{ 1/a_n, 1/(2 a_n), 1/(3 a_n), \dots \} \subset P_{n-1}$. This is disjoint from the indeterminacy locus since it does not contain point  $-1/a_n$ in $P_{n-1}$. \qed

A computation as in the proof of Lemma 3.1 shows that $k_{\cal X}$ maps $Q\leftrightarrow C_3$ according to:
$$\eqalign{k_{\cal X} \ :\  &Q \ni \xi_1 \mapsto [1:a_0-\xi_1:0] \in C_3, \cr &C_3 \ni [x_0:x_1:0] \mapsto -x_1/x_0 \in Q.} \eqno{(3.3)}$$

\proclaim Lemma 3.3.  If $a_0 \ne 2/m$ for all $m>0$ then the indeterminacy locus of $k_{\cal X}$ and the forward orbit of $C_4$ under the induced map $k_{\cal X}$ are disjoint. If $a_0 = 2/m$ for some $m>0$, we have $k_{\cal X}^{2m-1} C_4 = e_{01}$. 

\noindent{\it Proof.} Since the forward image of $C_4$ is $[1:-1+a_0 :0]\in C_3$, using (3.3) we have that $k_{\cal X}^{2m-1}C_4 = [1: m a_0 -1 :0] \in C_3$. Since the unique point of indeterminacy in $C_3$ is $e_{01}$, for $C_4$ to be mapped to a point of indeterminacy, $a_0$ must satisfy $m a_0-1 = 1$ for some  $m\ge0$.\qed

The following theorem comes directly from previous Lemmas.

\proclaim Theorem 3.4. Suppose that $n$ is even and  $a_0 \ne 2/m$ for all integers $m\ge0$.  Then no orbit of an exceptional curve contains a point of indeterminacy. 

Let us recall the Picard group $Pic({\cal X})$, which is the set of all divisors in ${\cal X}$, modulo linear equivalence, which means that $D_1\sim D_2$ if $D_1-D_2$ is the divisor of a rational function.  $Pic({\bf P}^2)$ is 1-dimensional and generated by the class of any line (hyperplane) $H$, and a basis of $Pic({\cal X})$ is given by the class of a general hyperplane $H_{\cal X}:=\pi^*H$, together with all of the blowup fibers $E_1, Q,P_1,\dots,P_{n-1}$.  If $r$ is a rational function on ${\cal X}$, then the pullback $k_{\cal X}^*r:=r\circ k_{\cal X}$ is just the composition.  To pull back a divisor, we just pull back its defining functions.  This gives the pullback map $k_{\cal X}^*:Pic({\cal X})\to Pic({\cal X})$.  Thus from (3.1-2) we see that the pullback of $E_1$ is $E_1$ and the pulling back of most of basis elements are trivial, that is $k_{\cal X}^* P_j =0$ for all $j=1, \dots, n-2$. 

Next we pull back $H_{\cal X}$.  Since $k$ has degree $2n+1$ we have $k^*H=(2n+1)H$ in $Pic({\bf P}^2)$.  Now we pull back by $\pi_{\cal X}^*$ to obtain:
$$(2n+1)H_{\cal X} = \pi_{\cal X}^*(2n+1) H = \pi_{\cal X}^* ( k^* H) .\eqno{(3.4)}$$
A line is given by  $\{h:= \alpha_0 x_0 + \alpha_1 x_1 + \alpha_2 x_2=0\}$, so $k^*H$ is the divisor defined by $h\circ k=\sum_j \alpha_j k_j$.  To write this divisor as a linear combination of basis elements $H_{\cal X}, E_1,Q,P_1,\dots, P_{n-1}$, we need to check the order of vanishing of $h\circ k$ at all of these sets.   Let us start with the coordinate system  $\pi_{\cal X}(\xi, v) = [v \xi:1:v]$  near $E_1$, defined in \S 2.  Using the expression for $k$ given in \S1 we see that $\alpha_0 k_0 + \alpha_1 k_1 + \alpha_2 k_2$ vanishes to order $n$ in $v$.  It follows that $ \pi_{\cal X}^* k^* H $ vanishes at $E_1$ with multiplicity $n$.  Similar computations for all other basis elements gives us $\pi_{\cal X}^* k^{*}  H=k_{\cal X}^*H_{\cal X} +nE_1+(n+1)Q+(n+1)\sum_j jP_j$. Combining with (3.4) we have
$$k_{\cal X}^*H_{\cal X} = (2n+1) H_{\cal X} -nE_1 -(n+1)Q-(n+1) \sum_{j=1}^{n-1}\, j P_j.\eqno{(3.5)}$$
Similarly,  we obtain:
$$\eqalign{k_{\cal X}^* \ :\ &Q \mapsto  H_{\cal X} -E_1-Q-P_1-2 P_2- \cdots - (n-1)P_{n-1} \cr &P_{n-1} \mapsto 2H_{\cal X}-E_1-Q-P_1-2 P_2- \cdots - (n-1)P_{n-1}.}\eqno{(3.6)}$$

\proclaim Theorem 1: ${\bf n=}$even. Suppose $F(z) = \sum_{j=1}^n a_j z^j$ is an even degree polynomial associated with $\jmath_F$. If $a_0 \ne 2/m$ for any positive integer $m$, then the degree complexity is the largest root of the quadratic polynomial $x^2-(n+1)x-1$.

\noindent{\it Proof.} Since $P_1, \dots, P_{n-2}$ are mapped to $0$ under the action on cohomology, it suffices to consider the action restricted to $H_{\cal X}, E_1, Q,$ and $P_{n-1}$. By (3.5,6) the matrix representation of $k_{{\cal X}}^* $, restricted to the ordered basis $\{H_{\cal X}, E_1,Q,P_{n-1}\}$, is
$$ \left( \matrix{ 2n+1 & 0&1&2\cr -n & 1& -1& -1\cr -n-1& 0 & -1&-1\cr -n^2+1& 0& -n+1&-n+1\cr} \right).$$
The characteristic polynomial  is $x(x-1)(x^2-(n+1) x-1)$.\qed

\bigskip
\noindent{\bf\S4. Mappings with $n=$ odd.}  Let us start with the space ${\cal X}$ from \S3.  When $n$ is odd, we see from (3.2) that the image of all exceptional lines of $k_{\cal X}$ coincide with a point of indeterminacy in $p_n\in P_{n-1}$.  Let $\pi_{\cal Y} : {\cal Y} \to \P^2$ be the complex manifold obtained by blowing up ${\cal X}$ at the point $p_n$, and let $P_n$ denote the exceptional fiber over $p_n$. In the $u_{n-1},\eta_{n-1}$ coordinate system, $p_n$ has coordinate $(0,1/a_n)_{n-1}$.  Thus, at $P_n$, we use the coordinate projection:
 $$\pi_n:{\cal Y}\ni (u, \eta)_n\to (u^{n}(u\eta+ 1/a_n), u^{n-1} (u\eta+1/a_n)) \in \C^2.$$ 
 Most computations in the previous section remain valid for $n$ odd.  Thus Lemma 3.3, (3.1) and (3.3) are still valid for the induced map $k_{\cal Y}: {\cal Y} \to {\cal Y}$. Under $k_{\cal Y}$ curves $C_1, C_2, P_1, \dots, P_{n-3}$ are still exceptional: 
$$\eqalign{& k_{\cal Y}\  :\  C_1, C_2, P_1 , \dots, P_{n-3} \mapsto -a_{n-1}/a_n^2 \in P_n\cr &k_{\cal Y}^{-1}\ :\ C_1, P_1, \dots, P_{n-3} \mapsto (a_{n-1}-(n-1)a_n)/a_n^2 \in P_n.}\eqno{(4.1)}$$
The blowup fibers $P_n$ and $P_{n-2}$ form a two cycle, $k_{\cal Y} : P_n \leftrightarrow P_{n-2}$ and  $P_{n-1}$ is mapped to itself as before. It follows that the points of indeterminacy for $k_{\cal Y}$ are $e_2, e_{01}$ and $ (a_{n-1}-(n-1)a_n)/a_n^2 \in P_n$.  For all $m \ge 0$, we have
$$k_{\cal Y}^{2m} \ : \ P_n \ni -a_{n-1}/a_n^2 \mapsto (2m(n-1)a_n - (4m+1) a_{n-1})/a_n^2 \in P_n\eqno(4.2)$$ 
As a consequence of (4.1) and (4.2) we have:
\proclaim Lemma 4.1.  If $n$ is odd, and if 
$$2a_{n-1} \ne (n-1) a_n,\eqno(4.3)$$ 
then the forward orbits of $C_1, C_2, P_1, \dots, P_{n-3}$ under $k_{\cal Y}$ do not contain any point of indeterminacy.

\medskip
\epsfysize=1.8in
\centerline{ \epsfbox{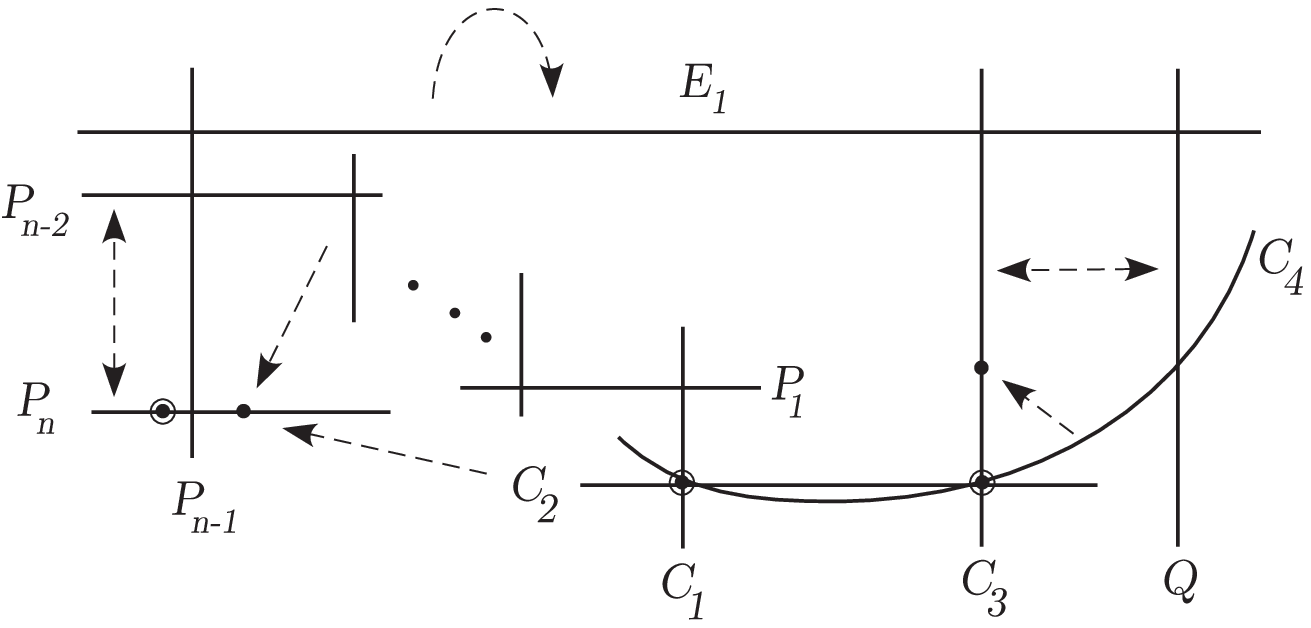}  }
\centerline{Figure 4.1.  The space ${\cal Y}$ and the action of $k_{\cal Y}$.}
\medskip

Combining Lemma 3.3 and 4.1 we have 
\proclaim Theorem 4.2. Suppose that $n$ is odd, $a_0 \ne 2/m$ for all  $m>0$, and $a_{n-1} \ne (n-1) a_n/2$.  Then the forward orbits of exceptional curves do not contain any points of indeterminacy. 
 
To determine $k_{\cal Y}$, we use the basis $\{H_{\cal Y},E_1,Q,P_1,\dots,P_n\}$ for $Pic({\cal Y})$.  Now the exceptional lines $C_1,C_2,P_1, \dots, P_{n-2}$ are mapped to $P_n$.  Let $\{C_1\}\in Pic({\cal Y})$ denote the class of the strict transform of $C_1$, i.e., the closure in ${\cal Y}$ of $\pi_{\cal Y}^{-1}(C_1-{\rm \ centers\ of\ blowup})$.  (The curve $C_2$ does not pass through any center of blowup, so with the same notation we have $\{C_2\}=H_{\cal Y}\in Pic({\cal Y})$.)   In order to write $\{C_1=(x_0=0)\}$ in terms of our basis, we note first that  $\pi_{\cal Y}^{-1} C_1 = C_4\cup E_1\cup Q\cup P_1 \cup \cdots \cup P_{n-1}$, i.e., the pullback function $x_0\circ\pi_{\cal Y}$ vanishes on all of these curves.  Thus we have to compute the multiplicities of vanishing.  At $P_{n-1}$, for instance, we consider the $(u_{n-1},\eta_{n-1})$ coordinate system defined in (2.4), and we see that $k_{\cal Y}^*x_0$ vanishes to order $n$ at $P_{n-1}=(u_{n-1}=0)$.   Similarly we can compute the multiplicities for $E_1, Q, P_1, \dots, P_{n-2}$ and $P_n$, so
$$H_{\cal Y}=\pi_{\cal Y}^* C_1 = \{C_1\}+ E_1+Q+2P_1+3P_2 + \cdots+ nP_{n-1}+nP_n.$$ 
It follows that 
$$k_{\cal Y}^* P_{n} = \{C_1\}+\{C_2\}+\sum_{j=1}^{n-2} \,  P_j= 2H_{\cal Y} -E_1-Q-\sum_{j=1}^{n-2} \, j P_j-n P_{n-1}-nP_n.$$ 
For the rest of basis entries we have 
$$\eqalign{k_{\cal Y}^* \ \ :\  & H_{\cal Y} \ \mapsto\  (2n+1)H_{\cal Y}-nE_1-(n+1)Q-(n+1)\sum_{j=1}^{n-1} j P_j -n^2 P_n\cr 
&Q\ \mapsto\ H_{\cal Y}-E_1-Q-P_1-2P_2-\cdots -(n-1)P_{n-1}-(n-1)P_n,\cr
&E_1 \mapsto E_1, \quad P_{n-2}\ \mapsto\ P_n,\quad {\rm and\ \ } P_{n-1}\ \mapsto\ P_{n-1}.\cr}$$

\proclaim Theorem 1: ${\bf n=}$odd.   If $a_0 \ne 2/m$ for all $m>0$, then the degree complexity is the largest root of the cubic polynomial $x^3-n x^2-(n+1)x-1$.

\noindent{\it Proof.}  The classes of the exceptional fibers $P_1, \cdots, P_{n-3}$ are all mapped to $0$, and exceptional fibers $E_1$ and $P_{n-1}$ are simply interchanged. It follows that to get the spectral radius of $k_{\cal Y}^*$ we only need to consider $4\times4$ matrix with ordered basis $\{H_{\cal X}, Q, P_{n-2},P_{n}\}$ and the spectral radius is given by the largest root of $x^3 -nx^2-(n+1)x-1$.\qed

\noindent{\bf\S5. Degree 3: a family of automorphisms.} Let us consider the 2 parameter family of maps $k= \jmath_F \circ \iota$ where $F(z) = a z^3+a z^2+b z+2$ with $a \ne0$. We consider the complex manifold $\pi_{\cal Z} : {\cal Z} \to \P^2$ obtained by blowing up $6$ points $e_2,e_{01},p_4,p_5,p_6,r$ in the complex manifold ${\cal Y}$ constructed in \S4.  As we construct the blowups, we will let $E_2,E_{01}, P_4, P_5,P_6$ and $R$ denote the exceptional fibers over $e_2,e_{01},p_4,p_5,p_6,$ and $r$ respectively.  Specifically, we blow up $e_2$ and $e_{01}$ and then:
$$\eqalign{&p_4:= -1/a \in P_3, \ \ p_5:= (2-b)/a \in P_4,\ \cr &p_6:=(2b-2-a)/a^2 \in P_5,\ \ {\rm and\ \ } r:=0 \in E_2 \cap \{x_1=0\}.}$$
We define the local coordinate system in a similar way we define local coordinates in \S 2. Using these local coordinates we can easily verify that  under the induced map $k_{\cal Z}$ we have
$$C_1 \to P_4 \to C_1, \  E_2 \to P_5 \to E_2, \  C_4 \to E_{01} \to C_4',\ {\rm and\ }C_2 \to P_6 \to R \to C_2'$$
and all mappings are dominant and holomorphic. 

%\medskip
%\epsfysize=1.9in
%\centerline{ \epsfbox{figures/spaceZv1.eps}  }
%\centerline{Figure 5.1.  Space ${\cal Y}$ with blowup fibers}
%\medskip

%
\medskip
\epsfysize=1.9in
\centerline{ \epsfbox{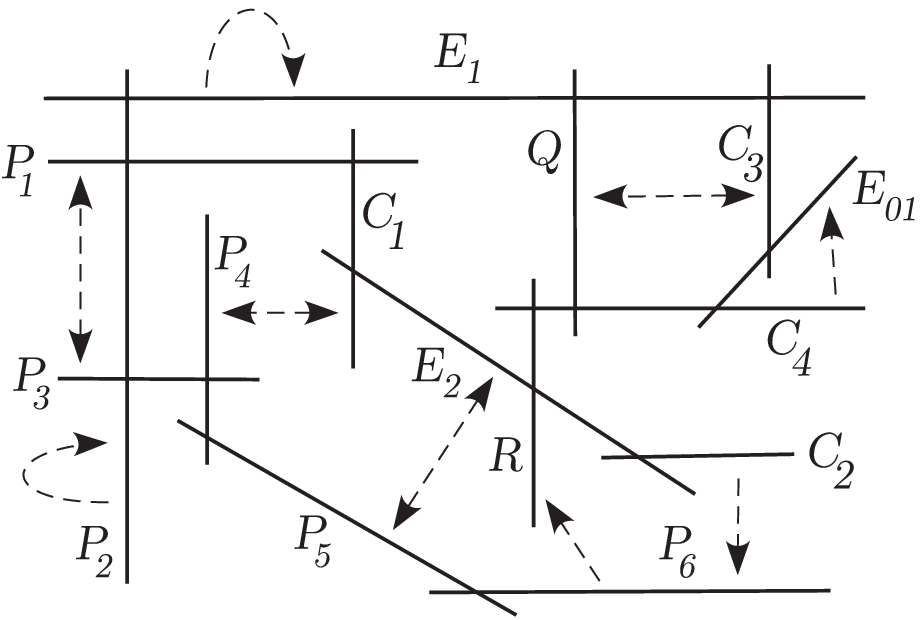}  }
\centerline{Figure 5.1.  The space ${\cal Z}$ and action of $f_{\cal Z}$.}
\medskip

\noindent For example, let us consider $E_2$.  We may use coordinates $w,\zeta$ which are mapped by $\pi_{E_2}:(x,\zeta) \to [w:w \zeta :1] \in \P^2$.  Thus  $E_2=(w=0)$ is given by $\zeta$-axis in this coordinate system and by considering $\lim_{w\to0} \pi_{P_5}^{-1}\circ k\circ \pi_{E_2}(w,\zeta)$ we find:
$$k_{\cal Z} \ :\  E_2 \ni \zeta \mapsto (2b-a-\zeta-1)/a^2 \in P_5.$$
The mapping among the exceptional fibers is shown in Figure 5.1.  What is not shown is that $R\to C_2'$ and  $E_{01}\to C_4'$

\proclaim Theorem 5.1. Suppose $F(z) = a z^3+a z^2+b z+2$ with $a \ne0$. Then the induced map $k_{\cal Z}$ is biholomorphic. 

\noindent{\it Proof.} Since $k_{\cal Z}$ and $k_{\cal Z}^{-1}$ have no exceptional hypersurface, indeterminacy locus for $k_{\cal Z}$ is empty. It follows that $k_{\cal Z}$ is an automorphism of ${\cal Z}$. \qed

Repeating the argument in previous two sections,  we have that $k_{\cal Z}^*$ acts on each basis element as follows :
$$\eqalign{&H_{\cal Z}\mapsto  7H_{\cal Z}-3E_1-4P_1-8P_2-9P_3 -10 P_4-10 P_5-10 P_6-3E_2-6R-4 Q - 4 E_{01},\cr &E_1 \mapsto E_1,\ \  P_1\mapsto P_3 \mapsto P_1,\ \ {\rm and\ \ }P_2 \mapsto P_2,\cr &P_4 \mapsto H_{\cal Z} -E_1-2 P_1-3P_2-3P_3-3P_4-3P_5-3P_6-E_2 - R - Q,\cr &P_5 \mapsto E_2, \ \  P_6 \mapsto H_{\cal Y} - E_2-R -E_{01},\ \  E_2 \mapsto P_5,\ \ \ {\rm and\ \ } E_{01} \mapsto P_6,\cr &Q \mapsto  H_{\cal Z} -E_1- P_1-2P_2-2P_3-2P_4-2P_5-2P_6 - Q-E_{01},\cr &E_{01} \mapsto  2H_{\cal Z} -E_1- P_1-2P_2-2P_3-2P_4-2P_5-2P_6-E_2 -2 R -2 Q-E_{01}.\cr  }$$

\proclaim Theorem 5.2.  Suppose $F(z) = a z^3+a z^2+b z+2$ with $a \ne0$. Then the degree of $k^n = k \circ \cdots \circ k$ grows quadratically, and $k$ is integrable. 

\noindent{\it Proof.}  All the eigenvalues of the characteristic polynomial of $k^*_{\cal Z}$ have modulus one.  The largest Jordan block in the matrix representation of $k_{\cal Z}^*$ is a $3 \times 3$  block corresponding to the eigenvalue $1$.  Thus the growth rate of the powers of the matrix is quadratic.  

Integrability follows from more general results: Gizatullin [G] showed that if the growth rate is quadratic, then there is an invariant fibration by elliptic curves.  In this case, we can give an explicit invariant.  If we 
define $\phi= \phi_1/\phi_2$ to be the quotient of the following two polynomials;
$$\eqalign{&\phi_1[x_0:x_1:x_2] = x_0^2x_2^2,\cr&\phi_2[x_0:x_1:x_2] = -2x_0^4+4x_0^3x_1-(2+a)x_0^2x_1^2+2a x_1x_2^2(x_0+x_2)-2b (x_0^3x_2- x_0^2x_1x_2),}$$  
then $\phi\circ k = \phi$. \qed

\centerline{\bf References}
{%\smrm 
\item{[A1]}    N. Abarenkova,  J-C. Angl\`es d'Auriac, S. Boukraa,
    S. Hassani and J-M. Maillard,  From Yang-Baxter equations to dynamical zeta functions for birational transformations.  Statistical physics on the eve of the 21st century, 436--490, Ser. Adv. Statist. Mech., 14, World Sci. Publishing, River Edge, NJ, 1999.

\item{[A2]}  N. Abarenkova,  J-C. Angl\`es d'Auriac, S. Boukraa,
    S. Hassani and J-M. Maillard,
    Rational dynamical zeta functions for
     birational transformations.
    Physica {A 264} (1999) pp. 264--293.  chao-dyn/9807014.

\item{[A3]} N. Abarenkova, J-C. Angl\`es d'Auriac,
   S.~Boukraa, S. Hassani and  J-M. Maillard,
   Topological entropy and complexity for discrete dynamical
   systems,     Phys.\ Lett.\ {A 262} (1999) pp. 44--49.  chao-dyn/9806026 .

\item{[A4]}  N. Abarenkova, J-C. Angl\`es d'Auriac, S. Boukraa and J-M. Maillard, Growth complexity spectrum of some discrete dynamical systems.  Phys. D 130 (1999), no. 1-2, 27--42.

\item{[A5]}  N. Abarenkova, J-C. Angl\`es d'Auriac, S. Boukraa and J-M. Maillard, Real topological entropy versus metric entropy for birational measure-preserving transformations.  Phys. D 144 (2000), no. 3-4, 387--433.

\item{[A6]}  N. Abarenkova, J-C. Angl\`es d'Auriac, S. Boukraa, S. Hassani and J-M. Maillard, Real Arnold complexity versus real topological entropy for birational transformations.  J. Phys. A 33 (2000), no. 8, 1465--1501.

\item{[A7]}  N. Abarenkova, J-C. Angl\`es d'Auriac, S. Boukraa, S. Hassani and J-M. Maillard, Topological entropy and Arnold complexity for two-dimensional mappings.  Phys. Lett. A 262 (1999), no. 1, 44--49.

\item{[A8]} N. Abarenkova,  J-C. Angl\`es d'Auriac, S. Boukraa,  and J-M. Maillard,
    Elliptic curves from finite order recursions or
   non-involutive permutations for discrete dynamical systems and
   lattice statistical mechanics,   The European  Physical Journal B. (1998), pp. 647--661.

\item{[BD1]} E. Bedford and J. Diller,  Real and complex dynamics of a family of birational maps of the plane: the golden mean subshift,  Amer. J. Math. 127 (2005), no. 3, 595--646.

\item{[BD2]}  E. Bedford and J. Diller,  Dynamics of a Two Parameter Family  of Plane  Birational Maps:  Maximal entropy,  J.\ of Geom. Analysis,  16 (2006), no. 3, 409--430. 

\item{[BD3]} E. Bedford and J. Diller,  Real dynamics of a family of plane birational maps: trapping regions and entropy zero, arXiv.math/0609113.

\item{[BK]} E. Bedford and KH Kim, Dynamics of rational surface automorphisms: Linear fractional recurrences.   arXiv:math/0611297

%\item{[BTV]} M. Bernardo, T. Truong and G. Rollet, The discrete Painlev\'e I equations: transcendental integrability and asymptotic solutions. J. Phys. A 34 (2001), no. 15, 3215--3252.

\item{[BV]} M. Bellon and C. Viallet, Algebraic entropy. Commun. Math. Phys. 204, 425--437 (1999)

%\item{[BHM1]}  S. Boukraa,   S. Hassani and  J-M. Maillard,  New integrable cases of a Cremona transformation: a finite-order orbits analysis, Physica A 240 (1997), 586--621.

\item{[BHM]}  S. Boukraa,   S. Hassani and  J-M. Maillard,  Product of involutions and fixed points, Alg.\ Rev.\ Nucl.\ Sci., Vol. 2 (1998), 1--16.

\item{[BM]} S.~Boukraa and  J-M. Maillard, Factorization properties of birational
     mappings,  Physica A {\bf 220} (1995),   pp. 403--470.
     
\item{[BMR1]}  S. Boukraa, J-M. Maillard, and G. Rollet, Almost integrable mappings,  Int.\ J. Mod. Phys. {B8} (1994), pp. 137--174.

\item{[BMR2]}  S. Boukraa, J-M. Maillard, and G. Rollet,  Integrable mappings and  polynomial growth,
   Physica A {\bf 209}, 162-222 (1994).

\item{[BMR3]}  S. Boukraa, J-M. Maillard, and G. Rollet,   Determinental identities on
     integrable mappings,  Int. J. Mod. Phys. {\bf B8} (1994), pp. 2157--2201.

\item{[DF]}  J. Diller and C. Favre,  Dynamics of bimeromorphic maps of surfaces, Amer. J.  Math., 123 (2001), 1135--1169.

%\item{[DS]} T-C. Dinh and N. Sibony,  Une borne sup\'erieure pour l'entropie topologique d'une application rationnelle.  preprint.

%\item{[D]}  R. Dujardin,  Laminar currents and birational dynamics.  preprint

\item{[FS]} J-E. Forn\ae ss and N. Sibony, Complex dynamics in higher dimension, II, { Modern Methods in Complex Analysis},  Ann.\ of Math.\ Studies, vol. 137, Princeton U. Press, 1995, p. 135--182.

\item{[G]} M. Gizatullin,  Rational G-surfaces, Izv. Akad. Nauk SSSR Ser. Mat. (44) 1980 110--144. 

%\item{[GH]} P. Griffiths and J. Harris, Principles of algebraic geometry. Wiley Classics Library, 1994.

%\item{[LM]}  D. Lind and B. Marcus,  An introduction to symbolic dynamics and coding. Cambridge University Press, Cambridge, 1995.

%\item{[R]} C. Robinson,  Dynamical systems. Stability, symbolic dynamics, and chaos. Second edition. Studies in Advanced Mathematics. CRC Press, Boca Raton, FL, 1999.

\item{[T1]}  T. Takenawa, A geometric approach to singularity confinement and algebraic entropy, J. Phys.\ A: Math.\ Gen.\ 34 (2001) L95--L102.

\item{[T2]} T. Takenawa, Discrete dynamical systems associated with root systems of indefinite type, Commun.\ Math.\  Phys., 224, 657--681 (2001).

\item{[T3]} T. Takenawa, Algebraic entropy and the space of initial values for discrete dynamical systems, J. Phys.\ A: Math.\ Gen.\ 34 (2001) 10533--10545.

\bigskip

\rightline{E. Bedford:  bedford@indiana.edu}

\rightline{T. Tuyen: ttuyen@indiana.edu}

\rightline {Department of Mathematics}

\rightline{Indiana University}

\rightline{Bloomington, IN 47405 USA}

\medskip\rightline{K. Kim: kim@math.fsu.edu}

\rightline {Department of Mathematics}

\rightline{Florida State University}

\rightline{Tallahassee, FL 32306 USA}

\medskip\rightline{N. Abarenkova:   nina@pdmi.ras.ru}

\rightline{Laboratory of Mathematical Problems of Physics}

\rightline{ Petersburg Department of Steklov Institute of Mathematics}

\rightline{27, Fontanka, 191023, St. Petersburg, Russia}

\medskip\rightline{J-M. Maillard:   maillard@lptmc.jussieu.fr}

\rightline{Lab. de Physique Th\'eorique et de la Mati\`ere Condens\'ee}

\rightline{Universit\'e de Paris 6, Tour 24}

\rightline{4\`eme \'etage, case 121, 4, Place Jussieu}

\rightline{75252 Paris Cedex 05, France}

}
\bye